\documentclass[a4paper,11pt]{article}
\usepackage[utf8]{inputenc}
\usepackage[round]{natbib}
\usepackage{url}

\usepackage{amsmath}
\usepackage{amsthm}
\usepackage{amssymb}

\usepackage{array}

\usepackage[american]{babel}
\usepackage{amsfonts}
\usepackage{algorithm}
\usepackage{algorithmic}
\usepackage{graphicx}
\usepackage{txfonts}
\usepackage{color}
\usepackage[normalem]{ulem}

\definecolor{gray}{rgb}{0.5,0.5,0.5}
\definecolor{red}{rgb}{1,0,0}

\usepackage
[
colorlinks=true,
linkcolor=blue,
anchorcolor=blue,
citecolor=blue,
urlcolor=blue,
]{hyperref}

\graphicspath{{./figs/}}

\newtheorem{theorem}{Theorem}
\newtheorem{corollary}[theorem]{Corollary}
\newtheorem{lemma}[theorem]{Lemma}
\newtheorem{definition}{Definition}
\newtheorem{proposition}[theorem]{Proposition}

\renewcommand{\O}[1]{\ensuremath{\mathcal{O}(#1)}}

\newcommand{\doi}[1]{DOI \href{http://dx.doi.org/#1}{\nolinkurl{#1}}}

\newcommand{\e}{{\mathbf{e}}}
\newcommand{\f}{{\mathbf{f}}}
\newcommand{\g}{{\mathbf{g}}}
\newcommand{\n}{{\mathbf{n}}}

\newcommand{\x}{{1}}
\newcommand{\y}{{2}}
\newcommand{\z}{{n}}

\newcommand{\ex}{{e}_\x}
\newcommand{\fx}{{f}_\x}
\newcommand{\ey}{{e}_\y}
\newcommand{\fy}{{f}_\y}
\newcommand{\ez}{{e}_\z}
\newcommand{\fz}{{f}_\z}

\newcommand{\gz}{{g}_\z}

\newcommand{\et}{{\e}_t}
\newcommand{\ft}{{\f}_t}
\newcommand{\gt}{{\g}_t}

\newcommand{\ed}{{e}_d}
\newcommand{\fd}{{f}_d}
\newcommand{\gd}{{g}_d}

\renewcommand{\k}{{\mathbf{k_e}}}

\newcommand{\W}{S}

\renewcommand{\S}{M}

\newcommand{\dk}{{\delta_\kappa}}

\renewcommand{\u}{{\mathbf{u}}}
\renewcommand{\v}{{\mathbf{v}}}

\newcommand{\vx}{{\mathbf{v}}_1}
\newcommand{\vy}{{\mathbf{v}}_2}

\newcommand{\kx}{{\kappa}_1}
\newcommand{\ky}{{\kappa}_2}
\newcommand{\kgx}{{\kappa}^g_1}
\newcommand{\kgy}{{\kappa}^g_2}

\newcommand{\te}{\theta_{\e}}

\newcommand{\K}{{\mathcal{K}}}
\newcommand{\dK}{{\mathcal{K}'}}

\newcommand{\fm}{\g}
\newcommand{\fp}{\f}

\newcommand{\Jm}{{\mathrm J}_\fm}
\newcommand{\Jp}{{\mathrm J}_\fp}
\newcommand{\J}{\Jp}

\newcommand{\Hp}{{H_p}}
\newcommand{\He}{{H_{\e}}}

\newcommand{\ke}{{k_{\e}}}

\newcommand{\Ae}{A_{\e}}

\newcommand{\Ap}{A_{p}}

\newcommand{\eg}{e.g.}
\newcommand{\ie}{i.e.}

\newcommand{\I}{\mathrm{I}}
\newcommand{\II}{\mathrm{I\!I}}

\theoremstyle{remark}

\numberwithin{remark}{section}
\newtheorem*{remark*}{Remark}

\begin{document}

\title{Uniform convergence of discrete curvatures\\ from nets of curvature lines}
\author{
Ulrich Bauer \thanks{Institute for Numerical and
Applied Mathematics,
University of G\"ottingen,
Lotzestr.~16--18,
37083~G\"ottingen, Germany. \href{http://ddg.math.uni-goettingen.de/}{\texttt{http://ddg.math.uni-goettingen.de/}}} \and Konrad Polthier
\thanks{
Department of Mathematics and Computer Science,
Freie Universit\"at Berlin,
Arnimallee 6,
14195~Berlin,
Germany. \href{http://geom.mi.fu-berlin.de/}{\texttt{http://geom.mi.fu-berlin.de/}}} \and Max Wardetzky \footnotemark[1] }
\date{December 11, 2009}

\maketitle

\begin{abstract}
We study discrete curvatures computed from nets of curvature
lines on a given smooth surface, and prove their uniform convergence
to smooth principal curvatures. We provide explicit error bounds, with
constants depending only on properties of the smooth limit surface and the shape regularity of the
discrete net.

\end{abstract}

\section{Introduction}




The 
field of {\em discrete differential geometry} has brought to light intriguing
discrete counterparts of classical differential geometric concepts as well as
efficient geometric algorithms~\citep[see, \eg,][]{DDGBook, SiggraphCourseNotes}.
One aspect of this theory is {\em convergence}: classical smooth
notions should arise in the limit of refinement. Recently,
several convergence results have been obtained for curvatures and differential
operators defined on polyhedral surfaces. Roughly, one may distinguish
three approaches: (i) polynomial
surface approximation~\citep[see, \eg,][]{MeekWalton00,CazalsPouget05}, (ii)
geometric measure theory~\citep[see,
\eg,][]{fu93convergence,MorvanCohen-Steiner06}, and (iii) finite element
analysis~\citep[see, \eg,][]{dziuk88laplace,wardetzky06convergence}. Among
these, (i) provides {\em pointwise} convergent curvatures for many, but not all,
discrete meshes. In contrast, (ii) and (iii) consider
generalizations of {\em integrated}, or total, curvatures and yield
convergence in the sense of measures or appropriate Sobolev norms,
respectively. 

Given the convergence of curvatures studied by approaches (ii) and (iii) in an {\em integrated} sense, it is natural to ask whether
these curvatures can be shown to
also converge in a {\em pointwise} manner. An affirmative answer can be obtained in
some special cases, such as polyhedral surfaces with vertices on the unit $2$-sphere~\citep{xu06sphere}. 
In general, however, the answer to this question
is \emph{negative}: it was observed in~\citet{xu05convergence} that for general
irregular polyhedral surfaces, there exist no {\em $k$-local} definitions of
discrete curvatures that are pointwise convergent. Here, by $k$-locality we mean
that the definition of curvatures associated with a vertex $p$ of a polyhedral surface only
depends on the $k$-star of $p$, \ie, those vertices that are connected to $p$ by a
path of at most $k$ edges. The concept of $k$-locality is motivated by the
smooth setting, where the definition of curvatures and differential operators only depends on local
properties of the underlying Riemannian manifold.

\begin{figure}[t]
\begin{center}
\includegraphics[height=3.5cm]{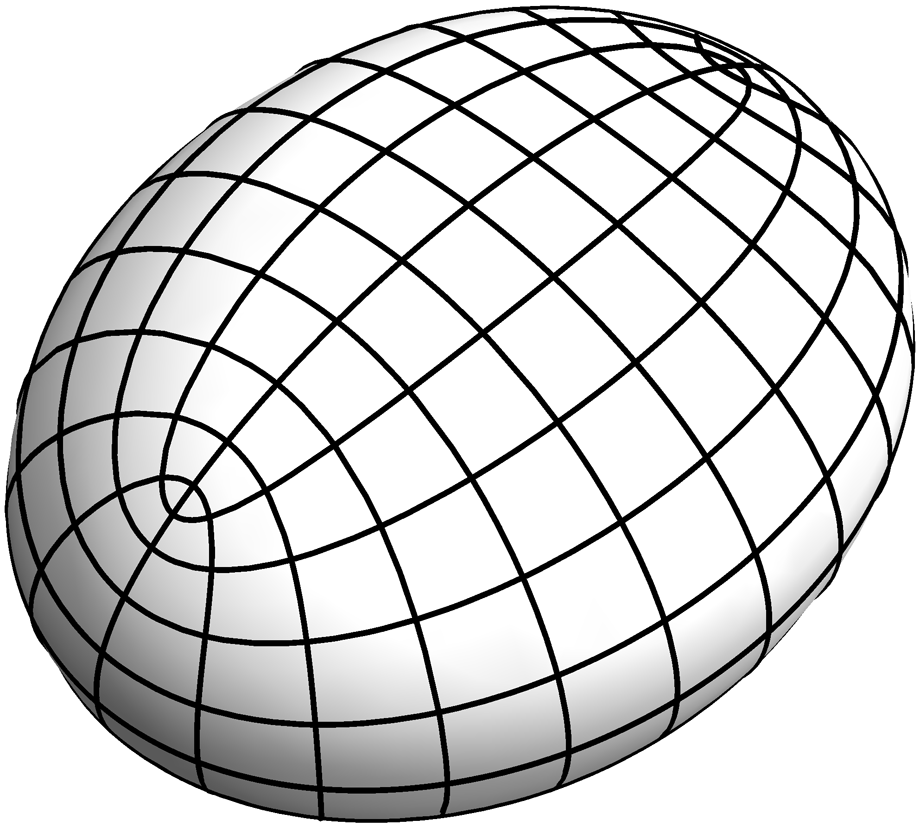}
\caption{A discrete net of curvature lines on an ellipsoid.}
\label{fig:ellipsoid}
\end{center}
\end{figure}

\paragraph*{Uniform convergence from nets of curvature lines}
We provide an affirmative answer to the above question of pointwise
convergence of curvatures for a {\em special class} of discrete meshes: discrete nets of
curvature lines on a given smooth surface $\S$ that is immersed into Euclidean space $\mathbb{E}^3$ (see
Figure~\ref{fig:ellipsoid}). To obtain approximations $(k_1, k_2)$ of principal curvatures $(\kappa_1, \kappa_2)$ on $\S$, we follow a three-step approach. We consider (a) local polyhedral approximations to nets of curvature lines, to which we apply (b) well-known $1$-local \emph{integrated} notions of discrete curvatures, such as those based on normal cycles~\citep[see, \eg,][]{MorvanCohen-Steiner06} or those based on the so-called \textit{cotangent formula}~\citep[see, \eg,][]{pinkall93minimal}, followed by (c) dividing the resulting integrated curvatures by appropriate area terms. The resulting pointwise curvature approximations $k_1,k_2:V\to \mathbb{R}$ are at first only defined on the vertex set $V$ of the underlying net. However, we may regard these curvature approximations as functions $k_i: \S\to \mathbb{R}$ by extending them from $V$ in a piecewise constant manner to the intrinsic Voronoi regions of the set $V$ on $\S$. 
Assuming this extension, we show:

\begin{theorem} \label{thm:convergence}
Let $\S$ be a smooth compact oriented surface without boundary\footnote{Surfaces with nonempty boundary can be treated with minor technical modifications.} immersed into $\mathbb{E}^3$. Consider a discrete net of curvature lines on $\S$ such that at each vertex the sampling condition~\eqref{eqn:samplingCondition} is satisfied. Let $\epsilon$ be an upper bound for the edge lengths of the net such that additionally the intrinsic $\epsilon$-balls around vertices cover all of $\S$. Then
\begin{align*}
\sup_{p\in\S}|k_i(p)-\kappa_i(p)| \leq C\epsilon\ , \quad i=1,2\ ,
\end{align*}
where $C$ depends only on properties of $\S$ and the shape regularity~\eqref{eqn:shapeRegularity} of the net of curvature lines.
\end{theorem}
A few remarks seem pertinent before proceeding:

\begin{list}
{\labelitemi}{\leftmargin=0em\itemindent=15pt}

\item \emph{Uniform} pointwise convergence of principal curvatures obtained by a $1$-local construction from nets of curvature lines is somewhat surprising since in general $1$-local polyhedral curvatures may not even converge in $L^2$, even if the mesh vertices reside on the smooth limit surface~\citep[see][]{wardetzky06convergence}. 

\item Our pointwise curvature approximations arise from dividing integrated ``Steiner-type'' curvatures by associated area terms. For example, the \emph{integrated} mean curvature of an interior edge of a polyhedral surface may be defined as the product between the length of that edge and the signed angle between the normals of its adjacent faces  -- a definition that arises from Steiner's view of considering offset surfaces. To obtain \emph{pointwise} curvature approximations from discrete integrated curvatures, we divide by so-called \emph{circumcentric areas}. While this approach is not new, see, \eg,~\cite{desbrun05dec}, our convergence result may be interpreted as a  justification of this construction \emph{provided} that the edges of a polyhedral surface well approximate the principal curvature directions of a smooth limit surface. We prove uniform lower and upper bounds for edge-based circumcentric areas that may be of interest in their own right.

\item   For our result to hold we require the explicit knowledge of positions of the vertices of a net of curvature lines on a smooth surface as well as the combinatorics of this net. More precisely, our curvature approximations at a vertex $p$ require the position of $p$ and the positions of its direct neighbors (with respect to the combinatorics of the net). Note that we do not require the knowledge of the entire net, though. (Given such an entire smooth net, it would be trivial to compute the principal curvatures at its vertices.) It would be desirable to drop from our approach the requirement of the exact knowledge of vertex positions of a smooth net of curvature lines. Here, one avenue for further study might be to consider discrete analogues of curvature line nets: so-called \emph{principle contact element nets}, see~\citet{Bobenko2008Discrete}.

\end{list}

\begin{figure}[t]
\begin{center}
\begin{tabular}{ccc}
\includegraphics[scale=1]{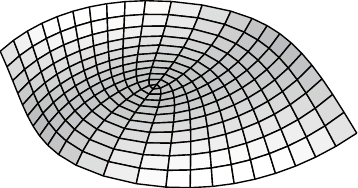}&
\lower2.5mm\hbox{\includegraphics[scale=1]{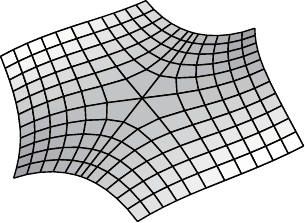}}&
\includegraphics[scale=1]{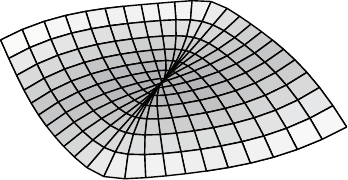}\\[9pt]
\emph{lemon}&\emph{star}&\emph{monstar}\\
\end{tabular}
\caption{The three generic patterns of curvature lines near an umbilical
point, called \emph{lemon, star}, and \emph{monstar} by~\citet[][]{berry77umbilic}.} \label{fig:umbilics}
\end{center}
\end{figure}

Our uniform convergence result given in Theorem~\ref{thm:convergence} is a consequence of a corresponding \emph{local} error estimate given in Theorem~\ref{thm:main}. This error estimate holds up to and including umbilical points, where singularities in the curvature line pattern arise. Using a refinement sequence for each of the three surfaces shown in Figure~\ref{fig:umbilics}, we observed numerically that although the shape regularity of the net may blow up near umbilics, linear convergence with respect to the maximum edge length $\epsilon$ remains valid in these cases. Our experiments also indicate that linear convergence is optimal.

\paragraph*{Alternative approaches}
An alternative point of departure for establishing pointwise convergence of
discrete curvatures is to give up $k$-locality and to
allow for $k\to\infty$ as the mesh refinement
increases. 
In fact, the above mentioned convergence results of (ii) and (iii) may be
interpreted in this way: by decreasing the diameter of the domains over which
discrete curvatures are integrated (measured), while simultaneously increasing
the mesh refinement inside these domains at a sufficiently fast rate, one
recovers classical pointwise notions of smooth curvatures in the limit.
In a similar fashion,~\citet{belkin08laplace} proposed a discrete Laplace operator based on the heat kernel. This operator converges in a pointwise manner if the kernel is scaled down while the mesh resolution is increased sufficiently fast relative to the scaling of the kernel. In contrast to
these works, which need to allow for $k\to\infty$ to establish pointwise
convergence, our result is obtained by working with the simplest and most
local definition: $k=1$.

\section{Discrete curvatures from nets of curvature
lines}\label{sec:discrete.curvatures}


In order to motivate our definition of discrete curvatures for nets of
curvature lines, we recall some important notions of
curvature for polygonal curves and polyhedral surfaces. For a similar
discussion, we refer to~\citet{sullivan08curves}.

\subsection{Discrete curvatures of ``Steiner-type''}\label{subsec:dicrete.curvatures}

\paragraph*{Integrated curvatures for polygonal curves}

Generalizations of classical smooth notions of curvature date at least back
to~\cite{steiner40parallel}, who considered parallel offsets of convex
hypersurfaces, relating {\em integrated} or {\em total} curvatures to changes
in length, area, and enclosed volume.
For example, for a convex curve $\gamma\subset \mathbb{E}^2$, one of Steiner's formulas
reads
\begin{align}\label{eqn:Steiner}
l(\gamma_\epsilon)=l(\gamma)+\epsilon \int_\gamma \kappa(s) ds ~,
\end{align}
where $l$ is the length functional, $\kappa$ denotes the curve's curvature, and $\gamma_\epsilon$ is the offset curve obtained by displacing $\gamma$ along its normals by some constant amount $\epsilon$.

\begin{figure}[t]
\begin{center}
\includegraphics[scale=1]{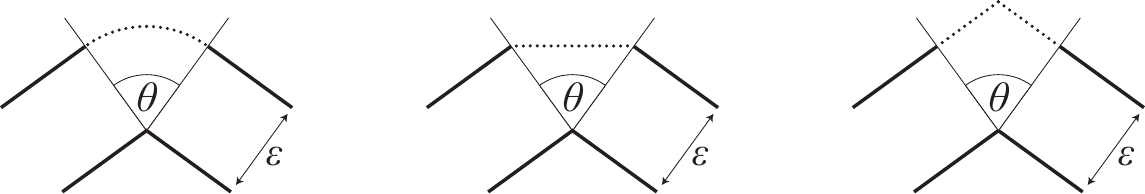}
\caption{Applying Steiner's formula~\eqref{eqn:Steiner} to the three depicted
definitions of offset curves for a given polygonal curve lead to the three
discrete curvatures in \eqref{eqn:polygonalVurvatures}.}
\label{fig:curvatureMeasure}
\end{center}
\end{figure}

Steiner's offset formula can be extended to the non-smooth and non-convex
case~\citep[][]{Federer59, Wintgen82,Zahle86}.
In particular, various notions for
curvatures of polygonal curves may be interpreted using Steiner's framework. Consider, \eg, 
\begin{align}\label{eqn:polygonalVurvatures}
     k_p \in \left\{\theta_p, ~~ 2\sin\frac{\theta_p}{2}, ~~ 2\tan\frac{\theta_p}{2}\right\} ~,
\end{align}
where $p$ denotes an inner vertex of a polygonal curve, and $\theta_p$ is the
turning angle between the two line segments incident to $p$. These notions arise by
applying~\eqref{eqn:Steiner} to the three different types of offsets
depicted in Figure~\ref{fig:curvatureMeasure}. Among these, the first notion
is the one considered by Steiner, the second corresponds to a finite element
discretization using piecewise linear functions, and the third also arises in the theory of discrete integrable
systems~\citep{bobenko99lagrange,hoffmann08hashimoto}.

\paragraph*{Integrated curvatures for polyhedral surfaces}
By a {\em polyhedral surface}, we mean a piecewise linear immersion of a compact simplicial surface into $\mathbb E^3$. Extending the
notions of discrete curvatures from polygonal curves to oriented polyhedral surfaces
leads to the following {\em edge-based} definitions of integrated normal
curvature:
\begin{align}\label{eqn:def.int.curv}
    \ke \in \left\{\te \|\e\|, ~~ 2\sin\frac{\te}{2} \|\e\|, ~~ 2\tan\frac{\te}{2}\|\e\|\right\}
    \ .
\end{align}
Here $\te\in(-\pi,\pi)$ is the signed angle between the normals of the two flat faces
incident to the edge $\e$. 
Notice that $\ke$ measures curvature \emph{orthogonal} to $\e$, since there is no curvature along $\e$ itself. Integrated {\em mean} curvature is accordingly defined as $\He=\frac{\ke}{2}$.

In the planar limit ($\te\to 0$), the definitions in~\eqref{eqn:def.int.curv} agree up to second
order in the angle variable. Therefore, as it turns out, it suffices for our purposes to
prove convergence of \emph{one} of these definitions in order to obtain
convergence for all of them. Convergence of the first definition
in~\eqref{eqn:def.int.curv} in the sense of {\em measures} was investigated
in~\citet{fu93convergence,MorvanCohen-Steiner06}.

For completeness, we remark that the above edge-based definitions give rise to {\em vertex-based} notions of
integrated {\em mean curvatures} by adding the mean curvatures over all edges emanating from a given
vertex $p$, \ie,
\begin{align}\label{eqn:def.int.curv.vertex}
    \Hp = \frac{1}{2}\sum_{\e\sim p} \He \ .
\end{align}
The factor $\frac12$ takes the meaning of distributing the normal
curvature of each edge equally among its two adjacent vertices.

\begin{figure}[t]
\begin{center}
\includegraphics[scale=1.3]{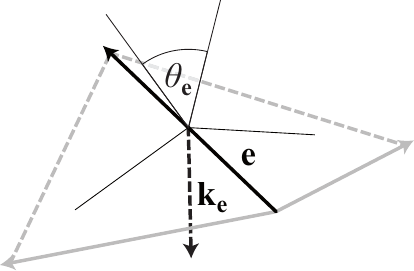}\hspace{3em}
\includegraphics[scale=1]{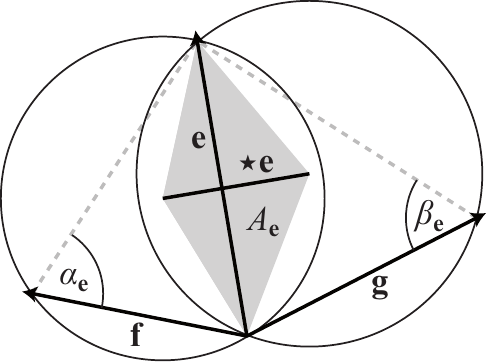}
\caption{Edge-based quantities. {\em Left}: dihedral angle $\theta_\e$ and discrete
curvature {vector}~$\k$. {\em Right}: dual edge ${\star\e}$ and circumcentric area $\Ae$.}
\label{fig:curvatureEdge}
\end{center}
\end{figure}

Finally, the scalar-valued definitions considered so far can be extended to corresponding
{\em vector-valued} notions. In the edge-based case, we obtain normal
curvature vectors $\k$ by multiplying $\ke$ with the angle-bisecting unit
normal vector at $\e$ (see Figure~\ref{fig:curvatureEdge}, left), and similarly for mean curvatures. Analogously to~\eqref{eqn:def.int.curv.vertex}, we then obtain vertex-based mean curvature vectors. We remark
that for $\ke = 2\sin\frac{\te}{2} \|\e\|$, the resulting mean curvature
vector coincides with the surface area gradient at $p$ when restricting to
piecewise linear surface variations~\citep[yielding the so-called \emph{cotangent
formula}, see~][]{pinkall93minimal}. Its convergence in the sense of Sobolev
norms was studied in~\citet{wardetzky06convergence}.

\paragraph*{From integrated to pointwise curvatures}

In order to obtain {\em pointwise curvatures}, we divide the above integrated
curvatures by corresponding area terms. Intuitively, these areas can be thought of as the domain of integration from which integrated curvatures were obtained.
Whether or not one obtains convergent curvatures depends on a careful choice of these areas. It turns out that for {\em triangulated}
polyhedral surfaces, one good choice are the so-called {\em circumcentric} areas, such as considered in~\citet{desbrun05dec}.
For each edge $\e$ we define
\begin{align}
\label{eqn:def.circumcentric}
    \Ae = \frac{1}{2} \mathrm{sgn}(\star e) \, \|\e\| \, \|{\star\e}\| \ ,
\end{align}
where $\|{\star\e}\|$ denotes the intrinsic length
of the circumcentric dual edge ${\star\e}$. This dual edge intrinsically
connects the circumcenters $C_1$ and $C_2$ of the two triangles $T_1$ and
$T_2$ incident to $\e$. Here, {\em intrinsic} means that one can think of
$T_1$ and $T_2$ as being unfolded onto the plane (see Figure~\ref{fig:curvatureEdge}, right). The
sign $\mathrm{sgn}(\star\e)$ is positive if along the direction of the ray from
$C_1$ through $C_2$, triangle $T_1$ lies before $T_2$, and negative otherwise. Note that $\mathrm{sgn}(\star\e) \leq 0$ (and therefore $\Ae \leq 0$) iff  $\alpha_\e +\,\beta_\e \geq \pi$, where $\alpha_\e$ and $\beta_\e$ are the angles opposite to $\e$ in the triangulation (see Figure~\ref{fig:curvatureEdge}, right).
Consequently, we require lower bounds that ensure positivity of circumcentric areas. For nets of curvature lines, we provide such bounds in Section~\ref{subsec:area}.

Similar to vertex-based integrated curvatures, we obtain vertex-based circumcentric areas from the edge-based case via
\begin{align}
\label{eqn:def.circumcentric.vertex}
    \Ap = \frac{1}{2} \sum_{\e\sim p} \Ae \ ,
\end{align}
where the sum is taken over all edges emanating from a given vertex $p$.
If all edges $\e$ incident to a vertex $p$ are intrinsically Delaunay (compare~\citet{BobenkoSpringborn05}), then $\Ap$ coincides with the intrinsic Voronoi area of $p$ and is therefore positive. 
However, as pointed out in~\citep{Dyer2009Circumcentric}, $\Ap$ might become negative in general, and a bound similar to the edge-based case is not possible. Therefore, we will not treat vertex-based pointwise curvatures based on $\Ap$.

\subsection{A local error estimate for discrete curvatures}\label{subsec:main}
In this section, we state our main local error estimate (Theorem~\ref{thm:main}), from which we derive our global uniform convergence result (Theorem~\ref{thm:convergence}).

Throughout we assume that $\S$ is a smooth compact oriented surface without boundary immersed into $\mathbb{E}^3$. By a {\em discrete net} on $\S$ we mean a cellular decomposition of $\S$ such that all attaching maps are homeomorphisms
and the intersection of any two cells is either empty or a single cell. As usual, we denote by $E$ the set of edges and by $V$ the set of vertices. We also assume that all edges are smoothly embedded. In a {\em discrete net of curvature lines} on $\S$, all edges are additionally required to be segments of curvature lines, non-umbilical vertices are required to have valence four, and umbilical vertices are required to have valence greater than two. In a completely umbilical region (such as $\mathbb{S}^2$), any net in the above sense serves as a net of curvature lines for our purposes. 
%




In order to be able to apply the concepts of discrete curvatures on polyhedral surfaces to nets of curvature lines, we require local polyhedral approximations of smooth curvature line nets.

\begin{figure}[t]
\begin{center}
\includegraphics[scale=1.2]{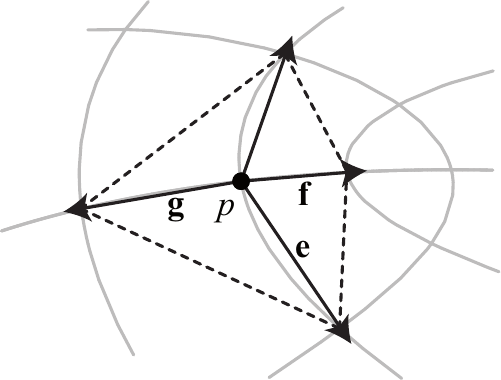}
\caption{The \emph{triangulated vertex star} at a vertex $p$ is a polyhedral surface approximating the (curved) discrete net of curvature lines at $p$.}
\label{fig:triangulatedVertexStar}
\end{center}
\end{figure}


\paragraph*{Local polyhedral approximation} %
In the sequel, the letters e, f, and g will be reserved for (curved) edges
in the edge set $E$, incident to a common vertex $p$, while corresponding bold face letters, $\e$, $\g$, and $\f$, will denote the {\em straight} edge vectors in
$\mathbb{E}^3$ obtained by connecting the endpoints of $e$, $f$, and $g$, respectively, by straight lines (see Figure~\ref{fig:triangulatedVertexStar}). We additionally assume that these edge vectors are \emph{oriented} such that they point away from $p$. For each disjoint edge pair $(e,f)$ incident to $p$
and contained in a common $2$-cell, we consider the flat triangle spanned by
$\e$ and $\f$. The union of these triangles forms
the triangulated {\em vertex star} of $p$, denoted by $st(p)$. 
Whenever we consider the triple $(\e,\f,\g)$, we will always assume that the pairs $(\e,\f)$ and $(\e,\g)$ span two triangles in $st(p)$, such that $\e$ becomes their common edge. Finally, as later justified by our sampling condition (\ref{eqn:samplingCondition}) and Corollary~\ref{cor:vertexStarProjection}, we may assume that $\n\cdot(\e\times\f)>0$ and
$\n\cdot(\e\times\g)<0$, where 
$\n$ denotes the normal of $\S$ at $p$. 



Each triangulated vertex star $st(p)$ thus yields the requisite \emph{local} polyhedral approximation\footnote{Observe that we do not require that our \emph{local} polyhedral approximations yield a consistent \emph{global} one.}, 
which forms the basis for our curvature approximations. 
As outlined in the previous section, our definition of pointwise curvatures relies on the division by certain circumcentric areas, which may become zero or negative in general. This motivates, for a given principal direction, to choose the associated edge vector with maximal circumcentric area. 

\begin{definition}[area maximizing edge]\label{def:areaMax}
Consider a vertex $p$ in a discrete net of curvature lines. If $p$ is umbilical, we call an edge vector $\e$  \emph{area maximizing} if it maximizes the circumcentric area among all edges emanating from $p$ in the local polyhedral approximation. If $p$ is non-umbilical, let $\vx$ be the principal direction canonically associated with an edge vector $\e$. We call $\e$ \emph{area maximizing} if it maximizes the circumcentric area among the two edge vectors associated with $\vx$. 
\end{definition}

We show in Section~\ref{subsec:area} that area maximizing edges always have circumcentric areas that are bounded away from zero.

\begin{definition}[principal curvature approximations]\label{def:discreteCurvature}
Consider a vertex $p$ in a discrete net of curvature lines and let $\e$ be an area maximizing edge (associated with a principal direction $\vx$ if $p$ is non-umbilical). Then 
\begin{align*}
k_2(p) :=\frac\ke{2\Ae}
\end{align*}
defines the \emph{principal curvature approximation} of $\ky(p)$, where $\ky$ refers to the principal curvature corresponding to $\vy$ if $p$ is non-umbilical and refers to the unique normal curvature if $p$ is umbilical. Here $\ke$ is one of the edge-based integrated polyhedral curvatures defined in~\eqref{eqn:def.int.curv}.
\end{definition}

The fact that in the above definition the edge vector $\e$ is associated with $\v_1$ while $\frac\ke{2\Ae}$ approximates $\kappa_2$ is not an oversight: $\ke$ measures curvature \emph{orthogonal} to $\e$.

\begin{figure}[t]
\begin{center}
\includegraphics[scale=.4]{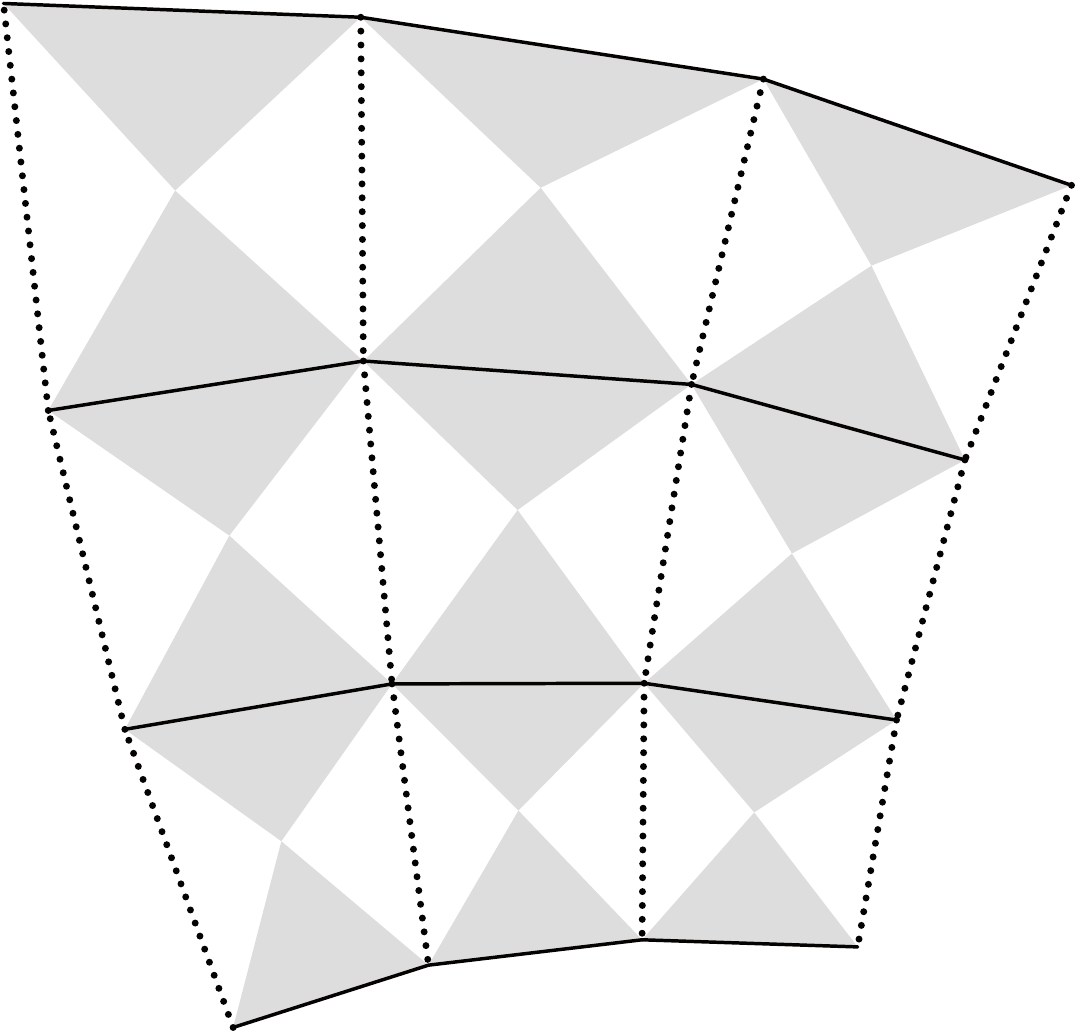}
\caption{A discrete net of curvature lines, with the two families of directions depicted as solid and dotted lines, respectively, together with the circumcentric areas (gray and white diamonds) per edge.}
\label{fig:tiling}
\end{center}
\end{figure}

The intuitive reason for dividing by \emph{twice} the circumcentric area in the above definition may (at least qualitatively) be explained as follows. Consider Figure~\ref{fig:tiling}, where one set of principal curvature directions, say those associated with $\v_1$, is depicted by solid lines, while the direction corresponding to $\v_2$ is represented by dotted ones. Likewise, the circumcentric areas corresponding to $\v_1$-directions are drawn as gray diamonds, while the circumcentric areas corresponding to $\v_2$-directions are represented as withe diamonds. Roughly, the gray diamonds cover only \emph{half} of the total surface. Hence, taking only the gray diamonds as regions of support for our (integrated) principal curvature approximations corresponding to $\kappa_1$, would mean to be roughly missing a factor of two. This motivates, for each edge along a given principal direction, to consider twice its circumcentric area as the domain of integration.

\paragraph*{Global constants}
For each $p\in\S$, let $\W(p)$ denote the shape
(or Weingarten) operator. Our estimates depend on both $\W$ and its covariant
derivative, $\nabla \W$. Accordingly, we define
\begin{align}\label{eqn:curvBound}
\K :=\mathop{\max_{p}}\| \W(p)\|_\text{op}%
\quad\text{and}\quad %
\dK := \mathop{\max_{p}}\left(\max_{{\|\v\|=1}}\|\nabla_\v S(p)\|_\text{op}\right) \ ,%
\end{align}
where $\|\cdot\|_\text{op}$ denotes the usual norm for linear operators.
Note that $\K$ is an upper bound for the normal curvatures of $\S$, whereas
$\dK$ provides an upper bound for directional derivatives $\nabla_{\v}(\kappa_i)$ of principal curvatures $\kappa_i$.

\paragraph*{Local constants} We also consider local constants -- shape regularity $\rho$ and maximum edge length $\epsilon$ -- that are specific for each vertex in the net of curvature lines. 
The reason for introducing local constants is that a high aspect ratio at one vertex should not affect the sampling condition (see below) at another vertex. In the following, we assume an arbitrary but fixed vertex $p$.

We let $\epsilon$ denote the largest {\em intrinsic} edge length over all (curved) edges $e\in E$ emanating from $p\in V$ (denoted by $e\sim p$),
\begin{align}\label{eqn:longestEdge}
\epsilon=\max_{e\sim p}l(e) \ .
\end{align}
Notice that the length of every edge vector emanating from $p$ in $st(p)$ is thus also bounded above by $\epsilon$.

Our estimates also depend on \emph{shape regularity}, or \emph{aspect
ratio}. We define $\rho\geq1$ to be the smallest number such that for \emph{all} pairs $(\e,\f)$ of edge vectors emanating from $p$ and forming a triangle in $st(p)$, one has
\begin{align}\label{eqn:shapeRegularity}
\frac\epsilon\rho\leq\|\e\|%
\quad\text{as well as}\quad
\frac{\|\e\|\,\|\f\|}{\|\e\times\f\|}\leq\rho %
\ .%
\end{align}
The former inequality implies $\frac1\rho\leq \frac{\|\e\|}{\|\f\|}\leq\rho$, while the latter means $\sin \angle(\e, \f) \geq \frac1\rho$.

\paragraph{Sampling condition}
In addition to the above definitions of maximum (local) edge length and (local) shape regularity, we assume the (local) sampling condition
\begin{align}\label{eqn:samplingCondition}
\epsilon \leq \frac{1}{16\K\rho^2} \ .
\end{align}
In some of our estimates, it will suffice to work with weaker sampling
conditions, such as $\epsilon \leq \frac{1}{2\K\rho}$ or $\epsilon \leq
\frac{1}{2\K}$, both of which are implied by \eqref{eqn:samplingCondition}.

\begin{theorem}[Local error estimate]
\label{thm:main} 
Let $\S$ be a smooth compact oriented surface without boundary immersed into $\mathbb{E}^3$. Consider a vertex $p\in V$ in a discrete net of curvature lines
on $\S$ and assume the
sampling condition~\eqref{eqn:samplingCondition}. Let $(k_1(p),k_2(p))$ denote the approximations of the smooth principal curvature $(\kappa_1(p),\kappa_2(p))$ as in Definition~\ref{def:discreteCurvature}. Then
\begin{align}\label{eqn:main.estimate}
\left|k_i(p) - \kappa_i(p) \right| \leq C\epsilon \quad i=1,2
\ .
\end{align}
The constant $C=C(\K,\dK, \rho)$ depends only on the curvature bounds~\eqref{eqn:curvBound} and the shape regularity~\eqref{eqn:shapeRegularity}. 
\end{theorem}


\begin{remark*}
Equivalent estimates can be obtained when replacing the scalar-valued
quantities in~\eqref{eqn:main.estimate} by
their corresponding vector-valued counterparts (for definitions, see
Section~\ref{sec:discrete.curvatures}). Proofs remain nearly identical.
\end{remark*}

Our global uniform convergence theorem stated in the introduction is a direct consequence of our local error estimate.  

\begin{proof}[Proof of Theorem~\ref{thm:convergence}]
Let $f:\S \to V$ be defined by mapping each $p\in \S$ to its nearest point in the vertex set $V$, i.e., $d_\S(p,f(p))=d_\S(p,V)$, where $d_\S$ denotes the intrinsic distance on $\S$.
We extend our curvature approximations (initially only defined on V) to functions $k_i: \S\to\mathbb{R}$ via $p \mapsto k_i(f(p))$. By the assumptions of Theorem~\ref{thm:convergence}, $\epsilon$ was \emph{globally} chosen such that $d_\S(p,f(p))\leq\epsilon$ for all $p\in M$. This implies, by connecting $p$ and $f(p)$ by a shortest geodesic arc, that
\begin{align*}
|\kappa_i(p)-\kappa_i(f(p))|\leq\dK\epsilon \ .
\end{align*}
Furthermore, $\epsilon$ was globally chosen such that it maximizes global edge length. Hence Theorem~\ref{thm:main} implies that there exists a constant $\tilde C$ such that 
\begin{align*}
 |k_i(f(p))-\kappa_i(f(p))|   \leq \tilde C\epsilon 
\end{align*}
for all $p\in \S$. Additionally, note that the constant $C$ in Theorem~\ref{thm:main}, besides depending on $\K$ and $\dK$, is monotonically increasing with respect to shape regularity. (This will become evident in the proof of Theorem~\ref{thm:main}.) Hence, $\tilde C$ depends only on $\K$, $\dK$, and the \emph{largest} local shape regularity constant $\rho$. This together with an application of the triangle inequality, 
\begin{align*}
|k_i(p)-\kappa_i(p)| = |k_i(f(p))-\kappa_i(p)|  \leq |k_i(f(p))-\kappa_i(f(p))| + |\kappa_i(f(p))-\kappa_i(p)| \ ,
\end{align*}
implies the claim.
\end{proof}

\section{Proof of local error estimate}\label{sec:proof}
The proof of Theorem~\ref{thm:main} proceeds in several steps. First, we
provide uniform lower and upper bounds for the edge-based circumcentric areas by which we divide integrated curvatures to obtain pointwise notions (Section~\ref{subsec:area}). In a
second step, we provide estimates for edge-based integrated curvatures by using their corresponding discrete curvature {\em
vector}. We establish that for each vertex $p$ in a net of curvature lines, the projection of these vectors onto the tangent plane
$T_p\S$ is negligible. Furthermore, we show that the remaining
normal component leads to the error estimate in Theorem~\ref{thm:main} up
to a certain error term (Section~\ref{subsec:Monge}). While for general meshes
the resulting error term cannot be controlled (and indeed causes failure of
pointwise convergence), we provide bounds for this error term for the
specific case of nets of curvature lines (Section~\ref{subsec:deltaKappa}).

%

\paragraph*{Basic assumptions} In order to avoid excessive repetition, we summarize our basic assumptions and notations. 
We write $p$ for a non-boundary vertex of a polyhedral surface, with vertex star denoted by $st(p)$. We assume that $\epsilon$ is an upper bound for the length of the (straight) edges emanating from $p$, and we let $\rho$ denote the shape regularity as defined in~\eqref{eqn:shapeRegularity}. If $st(p)$ arises from the local polyhedral approximation of a net of curvature lines, then $\epsilon$ is defined by~\eqref{eqn:longestEdge}, \ie,  as the maximum edge length of the curved edges emanating from $p$. Throughout, we assume the sampling condition~\eqref{eqn:samplingCondition}. As before, we denote curved edges by $e$, $f$, and $g$, and their corresponding straight edge vectors by $\e$, $\f$, and $\g$.


\subsection{Uniform bounds for circumcentric areas}\label{subsec:area}
In this section we prove upper and lower bounds for the circumcentric areas of area maximizing edges in the sense of Definition~\ref{def:areaMax}. 

\begin{proposition}\label{cor:largestEdgeIsFine}
Consider a vertex $p$ in a discrete net of curvature lines and let $\e$ be an area maximizing edge vector at $p$. Then our basic assumptions imply the existence of some $C>0$ such that
\begin{align*}
\frac {1}{C}\epsilon^2\leq \Ae \leq C\epsilon^2 \ ,
\end{align*}
where $C$ only depends on the shape regularity constant $\rho$.
\end{proposition}

Note that for non-umbilical vertices, this result implies the existence of an edge with positive circumcentric area for each of the two principal curvature directions. 

The remainder of this section is concerned with proving~Proposition~\ref{cor:largestEdgeIsFine}. First observe that
\begin{align*}
    \Ae = \frac{1}{4}(\cot \alpha_\e + \cot \beta_\e) \|\e\|^2 =
    \frac{\sin(\alpha_\e + \beta_\e)}{4\sin\alpha_\e\sin\beta_\e} \|\e\|^2\ ,
\end{align*}
where $\alpha_\e$ and $\beta_\e$ are the angles opposing $\e$ in the two
triangles meeting at $\e$, respectively. The requisite upper bound on $\Ae$ is relatively straightforward to obtain.

\begin{lemma}[upper bound]
\label{cor:circumcentricAreaUpperBound} Let $p$ be a vertex in a
discrete net of curvature lines on $\S$. Then our basic assumptions imply that each edge $\e$ emanating from $p$ satisfies 
\begin{align*}
\Ae \leq \rho^4\epsilon^2 \ .
\end{align*}
\end{lemma}
\begin{proof}
Clearly, we have
\begin{align*}
    \Ae \leq \frac{\|\e\|^2}{4\sin\alpha_\e\sin\beta_\e} \ .
\end{align*}
Let $\f$ be the edge vector emanating from $p$ such that $\alpha_\e$ belongs to the triangle formed by $\e$ and $\f$. Then the definition of shape regularity~\eqref{eqn:shapeRegularity}  implies
\begin{align}\label{eqn:angleBound}
\sin\alpha_\e =\frac{\|(\f-\e)\times\f\|}{\|(\f-\e)\|\,\|\f\|} \geq
\frac{\|\e\times\f\|}{(\rho+1)\|\e\|\,\|\f\|} \geq \frac{1}{(\rho+1)\rho} \geq
\frac{1}{2\rho^2} \ .
\end{align}
A similar estimate holds for $\beta_\e$. Hence, $\Ae \leq \rho^4 \| \e \|^2 \leq \rho^4 \epsilon^2$.
\end{proof}

Similarly, we obtain a lower bound for at least \emph{one} edge emanating from $p$.

\begin{lemma}[lower bound]
\label{cor:circumcentricAreaOneEdge} Let $p$ be a vertex in a
discrete net of curvature lines on $\S$. Then our basic assumptions imply that there exists an edge vector $\e$ emanating from $p$ such that
\begin{align*}
\Ae \geq \frac{1}{4\rho^3}\epsilon^2 \ .
\end{align*}
\end{lemma}
\begin{proof}
Let $\e$ be the shortest edge emanating from $p$, let $\f$ be the straight edge emanating from $p$ such that $\alpha_\e$ belongs to the triangle formed by $\e$ and $\f$, and define $\gamma:=\angle(\e,\f)$. Since $\|\e\|\leq \|\f\|$, it follows that $2\alpha_\e \leq (\pi-\gamma)$. Since in particular $0<\alpha_\e < \frac\pi2$, it follows that
$\cot \alpha_\e >\frac\pi2-\alpha_\e$. By~\eqref{eqn:shapeRegularity} we have $\sin\gamma \geq \frac1\rho$, and hence
\begin{align*}
\cot \alpha_\e >\frac\pi2-\alpha_\e
\geq \frac{\gamma}{2}
\geq \frac{\sin\gamma}{2}
\geq \frac{1}{2\rho}
~.
\end{align*}
Applying similar arguments, we obtain $\cot \beta_\e > \frac{1}{2\rho}$. Together, this yields
\begin{align*}
\Ae &= \frac{1}{4}(\cot \alpha_\e + \cot \beta_\e) \|\e\|^2
\geq \frac{1}{4\rho} \|\e\|^2
\geq \frac{1}{4\rho^3} \epsilon^2
~. \qedhere
\end{align*}
\end{proof}

\paragraph*{Lower bounds for vertices of valence four} The above lower bound on the circumcentric area of at least one edge suffices for umbilical vertices. However, it does not suffice at non-umbilical ones, since we require lower bounds for edges associated with \emph{each} of the two principal directions. We achieve this by showing that for each vertex of valence four, there are at least three edges that satisfy the required lower bound. Hence, for each vertex in a net of curvature lines, we have at least one good edge per principal direction. 

We note that some of the following results are also valid for vertices of valence different from four, such as, in particular, Corollary~\ref{cor:vertexStarProjection}.

As before, we let $\alpha_\e$ and $\beta_\e$ denote the angles opposing the straight edge $\e$ in the two
triangles meeting at $\e$, respectively. If $\alpha_\e \geq \delta$, $\beta_\e \geq \delta$, and $\alpha_\e+\beta_\e \leq\pi-\delta$ for some $\delta > 0$, then
\begin{align}\label{eqn:circumcentric.area}
    \Ae \geq \frac{\sin\delta \|\e\|^2}{4}  \ ,
\end{align}
which provides a useful lower bound if $\delta$ can be bounded away from zero.
Accordingly, we introduce the
notion of $\delta$-Delaunay edges, a nomenclature that is borrowed from the
classical case of {\em Delaunay triangulations} (corresponding to $\delta =
0$).

\begin{definition}[$\delta$-Delaunay]
Let $\alpha_\e$ and $\beta_\e$ be the angles opposing an edge $\e$ in the two
triangles meeting at $\e$, respectively. Then $\e$ is called $\delta$-Delaunay if there exists $\delta\geq 0$ such
that $\alpha_\e \geq \delta$, $\beta_\e \geq \delta$, and $\alpha_\e+\beta_\e
\leq\pi-\delta$. 
\end{definition}

Assume for a moment that $p$ has valence four and that $st(p)$ is planar. Assume further that all of the eight angles opposing the four edges emanating from $p$ are bounded from below by $\delta$. Then it is straightforward to verify that at least three among the four edges incident to $p$ are $\delta$-Delaunay. 
In general, however, $st(p)$ is not planar, and we have to account for Gaussian curvature. As usual, we define discrete Gauss curvature at a vertex $p$ of a polyhedral surface as the angle
defect, \ie, by $K_p=2\pi -\sum_i\gamma_i$, where $\gamma_i$ are the
intrinsic angles meeting at $p$. We obtain:

\begin{lemma}\label{lemma:Delaunay}
Let $p$ be a non-boundary vertex of valence four on a triangulated polyhedral surface, and let $\alpha_i, \beta_i$ denote the pairs of
angles opposing the four edges emanating from $p$. Assume that $\alpha_i
\geq\delta$ and $\beta_i\geq\delta$ for all $i=1,\dots,4$, as well as
$K_p<2\delta$. Then at least three among the four edges meeting at $p$ are
$\delta$-Delaunay.
\end{lemma}
\begin{proof}
Using $K_p = -2\pi+\sum_{i=1}^4 (\alpha_i+\beta_i)$, the result follows
from a straightforward calculation.
\end{proof}
%

We now show that our basic assumptions imply the assumptions of Lemma~\ref{lemma:Delaunay} when setting $\delta:=\frac{1}{2\rho^2}$. To see this, first observe
that with the above notations, $\alpha_\e\geq \sin\alpha_\e \geq \delta$ by~\eqref{eqn:angleBound}, and analogously for $\beta_\e$. It remains to check the condition  $K_p<2\delta$ for the discrete Gauss curvature. The requisite bound will be established in Lemmas~\ref{lemma:Gauss-bound} and~\ref{lemma:normalBound}. The resulting consequence for the existence of $\delta$-Delaunay edges is summarized in Lemma~\ref{lemma:delta.Delaunay}. Finally, Lemma~\ref{cor:circumcentricAreaEdge} establishes the lower bound for $\Ae$.

\begin{lemma}\label{lemma:Gauss-bound}
Let $p$ be a non-boundary vertex of an oriented polyhedral surface. Assume that the normals of the triangles incident to $p$ make an angle no
greater than $\phi\in [0,\pi)$ with some fixed
direction in $\mathbb{E}^3$. 
Then the discrete Gauss
curvature associated with $p$ satisfies $K_p\leq 2\pi (1-\cos\phi)$.
\end{lemma}

\begin{proof}
The lemma is a consequence of the Gauss-Bonnet theorem. Let $\n_1, \n_2, \dots, \n_l$ denote the unit normals of the triangles $T_1, T_2,\dots,T_l$ incident to $p$, ordered according to the orientation of  the polyhedral surface. Each $\n_i$ represents a point on
$\mathbb{S}^2$. Connecting consecutive pairs $(\n_i, \n_{i+1})$ 
by geodesic arcs on $\mathbb{S}^2$ yields a spherical polygon $P$ (possibly with intersecting edges). For each $i$, the exterior angle of $P$ at $\n_i$ is then equal to the interior angle of the Euclidean triangle $T_i$ at $p$. In particular, the sum $\Sigma$ of the exterior angles of $P$ satisfies $\Sigma=2\pi-K_p$. 

Moreover, by assumption we have $\phi\in [0,\pi)$, so all $\n_i$ lie on the same hemisphere. We can hence consider the spherical convex hull $\tilde P$ of $\n_1, \n_2, \dots, \n_l$, \ie, the smallest spherical polygon that contains all $\n_i$ and that is convex with respect to shortest geodesic arcs. Let $\tilde\Sigma$ denote the sum of exterior angles of $\tilde P$. It is easy to verify that $\tilde\Sigma\leq\Sigma$.  Hence $K_p \leq 2\pi -\tilde\Sigma=area(\tilde P)$, where the last equality follows from the Gauss-Bonnet theorem. Since $\phi\in [0,\pi)$, the polygon $\tilde P$ is contained in a geodesic disk of radius $\phi$, the area of which is $2\pi (1-\cos\phi)$. This proves the claim.
\end{proof}

In order to make use of the previous lemma, we seek a bound on the angle
$\phi$ between the surface normal at $p \in \S$ and the normals to the
triangles incident to $p$. Note that our basic assumptions, and in particular our sampling condition, imply $\epsilon \leq\frac{1}{2\K}$, where $\K$ is our curvature bound. Hence we can infer from~\citet[][Section~3, Corollary
1]{MorvanThibert04}:

\begin{lemma}\label{lemma:normalBound}
Let $p$ be a vertex (umbilical or not) in a net of curvature lines on $M$. Given our basic assumptions, $st(p)$ can be oriented such that the maximum angle $\phi\in [0,\pi)$
between the surface normal at $p\in \S$ and the normals to the triangles in $st(p)$ satisfies $\sin\phi \leq \left(4\rho+2\right)\K\epsilon$,
where $\rho$ is the shape regularity. In particular, $\sin\phi \leq \frac38$.
\end{lemma}

\begin{corollary}\label{cor:vertexStarProjection}
Under the assumptions of Lemma~\ref{lemma:normalBound}, $st(p)$ can be oriented such that the orthogonal projection of $st(p)$ onto $T_pM$ is injective and orientation-preserving.
\end{corollary}

Note that our assumptions are slightly different from those used
in~\citet{MorvanThibert04}; in fact, our assumptions are stricter. While their
sampling condition bounds the maximal extrinsic distance between two
neighboring vertices, we consider the intrinsic length on the smooth surface,
which is always larger. Additionally,~\citet{MorvanThibert04} require that the
distance between the discrete and the smooth surface is less than the {\em
reach} of the smooth one. This requirement is implicitly fulfilled locally by
the sampling condition $\epsilon \leq\frac{1}{2\K}$, since the reach of a surface patch formed by an intrinsic $\epsilon$-disk around $p$ is nothing but the minimal radius of curvature of that surface patch .

\begin{lemma}\label{lemma:delta.Delaunay}
In addition to the assumptions of Lemma~\ref{lemma:normalBound}, assume that $p$ is of valence four. Let $\delta:=\frac{1}{2\rho^2}$, where $\rho$ is the shape regularity. Then
at least three edges incident to~$p$ are $\delta$-Delaunay.
\end{lemma}
\begin{proof}
Using that $(1-\cos\phi)\leq \sin^2\phi$ for $\phi\in
[0,\frac{\pi}{2}]$, Lemmas~\ref{lemma:Gauss-bound}
and~\ref{lemma:normalBound} show that the discrete Gauss curvature at $p$
satisfies
\begin{align*}
K_p \leq 2\pi\left(4\rho+2\right)^2\K^2\epsilon^2 \leq (16 \K \rho \epsilon)^2
\ .
\end{align*}
Setting $\delta=\frac{1}{2\rho^2}$, the sampling
condition~\eqref{eqn:samplingCondition} implies $K_p \leq 2\delta$.
Moreover,~\eqref{eqn:angleBound} implies $\alpha_i \geq\delta$ and
$\beta_i\geq\delta$ for all pairs of angles $\alpha_i, \beta_i$ opposing the
four straight edges emanating from $p$. Finally, Lemma~\ref{lemma:Delaunay} implies the claim.
\end{proof}

\begin{lemma}[lower bound for valence four]
\label{cor:circumcentricAreaEdge} Under the assumptions of Lemma~\ref{lemma:delta.Delaunay}, there exist at least three straight edges among the four edges emanating from $p$, such that 
\begin{align*}
\Ae \geq \frac{1}{16\rho^4}\epsilon^2 \ ,
\end{align*}
for each edge $\e$ among these three. 
\end{lemma}
\begin{proof}
Observe that~\eqref{eqn:circumcentric.area} implies $4 \Ae \geq \sin\delta \|\e\|^2$ for all $\delta$-Delaunay edges. Applying Lemma~\ref{lemma:delta.Delaunay} and using $\sin\left(\frac{1}{2\rho^2}\right)\geq\frac{1}{4\rho^2}$ gives
\begin{align*}
\Ae &\geq \frac{1}{4}\sin\left(\frac{1}{2\rho^2}\right) \|\e\|^2
\geq \frac{\|\e\|^2}{16\rho^2} 
\geq \frac{1}{16\rho^4}\epsilon^2 \ . \qedhere
\end{align*}
\end{proof}

\begin{proof}[Proof of Proposition~\ref{cor:largestEdgeIsFine}]
Lemma~\ref{cor:circumcentricAreaUpperBound} provides an upper bound for both umbilical and non-umbilical vertices. Lemma~\ref{cor:circumcentricAreaOneEdge} provides the requisite lower bound for umbilical vertices. Finally, Lemma~\ref{cor:circumcentricAreaEdge} provides the lower bound for non-umbilical ones, since it implies the existence of at least one edge per principal direction with circumcentric area bounded from below.
\end{proof}

\subsection{Estimates for discrete integrated curvatures}
\label{subsec:Monge}

In this section, we establish a bound for the difference between the edge-based \emph{integrated} curvatures $\ke$ and (an appropriately scaled version of) the smooth principal curvatures of $\S$. Specifically, for a given edge $\e$ in our polyhedral approximation, we work with the discrete integrated curvature $\ke = 2\sin\frac{\te}{2} \|\e\|$ introduced in Section~\ref{subsec:dicrete.curvatures}.

\paragraph*{Darboux frames} %
For our purposes, it turns out to be useful to express vectors in frames that are locally adapted to the geometry of the surface $\S$. Specifically, a \emph{Darboux frame} at a non-umbilical point $p\in\S$ is an adapted frame given by $(\vx,\vy,\n)$, where $\vx$ and
$\vy$ are (normalized) principal directions of $\S$ at $p$, and $\n=\vx\times\vy$
is the surface normal (induced by the orientation of $\S$). For umbilical points, any adapted (\ie, $\vx,\vy\in T_p\S$) and orthonormal  (\ie, $\|\vx\|=\|\vy\|=1$ and $\n = \vx \times \vy$) frame may be considered a Darboux frame. 
Throughout, we employ the notation
\begin{align*}
\e = (\ex,\ey,\ez) 
\end{align*}
to represent a vector $\e$ in the coordinates given by a Darboux frame.
In the sequel, we assume a fixed Darboux frame at every vertex $p$ of our discrete net of curvature lines. If $p$ is non-umbilical then each edge vector emanating from $p$ is \emph{canonically associated} to exactly one of the principal directions $\vx$ or $\vy$. If $p$ is umbilical, we additionally require an explicit association of each edge vector with one of either $\vx$ or $\vy$.
In order to state the main result of this section, we require the notion of \emph{tangential deviation} of an edge vector with respect to a Darboux frame.

%

\begin{figure}[t]
\begin{center}
\includegraphics[scale=1.2]{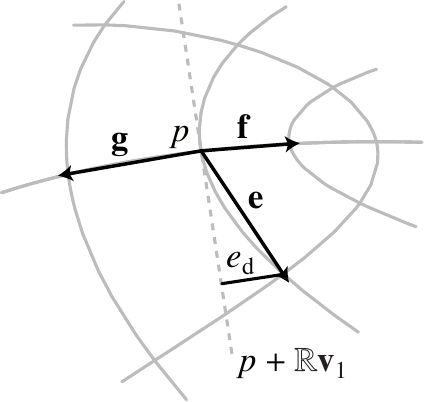}
\caption{The \emph{tangential deviation}~$\ed$ of a straight edge $\e$.}
\label{fig:tangentialDeviation}
\end{center}
\end{figure}

\begin{definition}[tangential deviation]\label{def:tangentialDeviation}
Let $p$ be a vertex in a discrete net of curvature lines on $\S$ and let $(\vx, \vy, \n)$ be the Darboux frame at $p$. If an edge vector $\e$ is associated with the principal direction~$\vx$, we call $\ed =  |\ey|$ its \emph{tangential deviation}, see Figure~\ref{fig:tangentialDeviation}. Likewise, if $\e$ is associated with $\vy$, we let $\ed =  |\ex|$. 
\end{definition}

The following proposition summarizes the main result of this section.


\begin{proposition}
\label{lemma:edgeCurvatureBound}
Using the notations of Definition~\ref{def:tangentialDeviation}, let the edge vector $\e$ emanating from vertex $p$ be associated with $\vx$, and let its directly neighboring edges $\f$ and $\g$ be associated with $\vy$. Then our basic assumptions imply that
\begin{align*}
\left| \ke - 2A_\e \ky \right|
\leq C\left(\left| \dk(\ed+\fd+\gd) \right| \epsilon + \epsilon^3 \right)  \ ,
\end{align*}
where $\kappa_2$ denotes the principal curvature of $\S$ in direction $\v_2$, $\dk = \ky-\kx$ (with $\dk=0$ in the umbilical case), and $C= C(\K, \dK, \rho)$.
\end{proposition}

\begin{remark*}
We point out that the estimates in this section are not entirely specific for nets of curvature lines. In fact, they hold in the more general setting of arbitrary smooth nets embedded in $\S$, as long as we require bounded shape regularity and the sampling condition $\epsilon\leq\frac{1}{2\K}$, where $\epsilon$ is the maximum intrinsic edge length of the net and $\K$ denotes our usual curvature bound. As a consequence, the estimates of the current section do not suffice to establish our main error estimate. Indeed, to obtain uniform convergence, we require that the term $\dk(\ed+\fd+\gd)$ appearing in Proposition~\ref{lemma:edgeCurvatureBound} is of order $\epsilon^2$. However, this is false for general nets (and causes failure of uniform convergence), but is true for nets of curvature lines as we will show in Section~\ref{subsec:deltaKappa}. For clarity's sake, though, we have decided to restrict the discussion of the current section to nets of curvature lines and rely on our (rather strong) \emph{basic assumptions} set forth in the beginning of Section~\ref{sec:proof}.
\end{remark*}

Proposition~\ref{lemma:edgeCurvatureBound} is proven in several steps. We commence by estimating the normal component $\ez$  of the straight edge vector $\e$ (Lemma~\ref{lemma:3rdOrderBound} and Corollary~\ref{edge normal component}). Then we switch to the edge-based curvature \emph{vector} $\k$ corresponding to $\ke$. We show that its tangential component is negligible (Lemma~\ref{tangentialBound}). The final proof of Proposition~\ref{lemma:edgeCurvatureBound} is given at the end of this section, where we show that the normal component of $\k$ yields the desired estimate.

\begin{lemma}\label{lemma:3rdOrderBound}
Let $p$ be a vertex of a discrete net of curvature lines on $\S$. Consider an edge $e\in E$ emanating from $p$ with corresponding straight edge vector $\e$. Writing $\e=(\ex,\ey,\ez)$ with respect to a Darboux frame centered
at $p$, our basic assumptions imply that
\begin{align*}
|\ez|\leq \K \e^2 \leq \K \epsilon^2\ .
\end{align*}
Furthermore, if $P(x,y)=\frac{\kx}{2}x^2+\frac{\ky}{2}y^2$ denotes the
osculating paraboloid, then
\begin{align*}
|\ez-P(\ex,\ey)|\leq C \epsilon^3 \ ,
\end{align*}
where $C$ only depends on our global curvature bound $\K$ and the bound on curvature derivatives $\dK$.
\end{lemma}

\begin{proof}
Using a Darboux frame at $p$, the surface $\S$ can locally be parameterized by a height function $h(x,y)$ over the tangent plane $T_p\S$. In the coordinates of the Darboux frame, we have $\e=(\ex,\ey,\ez)$, with $\ez=h(\ex,\ey)$. We let $d= \|(\ex,\ey)\|$, where throughout this proof $\|\cdot\|$ denotes the Euclidean norm in the parameter domain. Furthermore, we consider the \emph{constant} vector field $\v(x,y)=(\ex,\ey)/d$ in the parameter domain. 

Let $D_\v^{i} h$ denote the $i$th iteration of the directional derivative of $h$ along $\v$, \ie, $D_\v^1 h = D_\v h$ and $D^i_\v = D_\v (D^{i-1}_\v h)$, and let $D^i h$ denote the $i$th total derivative with respect to the standard Euclidean metric in the parameter domain. Observe that $h(0) = Dh(0)=0$, and hence
\begin{align}\label{eqn:intH}
\ez=h(\ex,\ey)&= \int_0^d\int_0^t D_\v^2h(\tau \v) \, d\tau \, dt \ .
\end{align}
Consequently, in order to prove the first part of the lemma, we seek an upper bound on $|D_\v^2 h|$ in terms of $\K$. 

Let $\W$ denote the shape operator and let $\I(\cdot,\cdot)$ and $\II(\cdot,\cdot)=\I(\W\cdot,\cdot)$ be
the first and second fundamental form of $\S$, respectively, with respect to the local parameterization induced by $h$. From
\begin{align} \label{eqn:I}
\I(\u,\v) &= \u^T (\mathrm{Id} + Dh \, Dh^T) \v 
\end{align}
and
\begin{align}
\label{eqn:kv}
\kappa_\v&=\frac{\II(\v,\v)}{\I(\v,\v)}=\frac{\I(\W\v,\v)}{\I(\v,\v)}=\frac{D^2_\v
h}{\sqrt{1+\|Dh\|^2} \I(\v,\v) } \ ,
\end{align}
we obtain the estimate 
\begin{align}
\label{eqn:d2vEstimate}
|D^2_\v h| = |\kappa_\v| \sqrt{1+\|Dh\|^2} \I(\v,\v) 
\leq  \K (1+\|Dh\|^2)^\frac{3}{2} \ .
\end{align}
In order to bound $\|Dh\|$, first observe that our sampling condition implies $\epsilon\leq\frac{1}{2\K}$. This, in turn, provides a bound on the (positive) angle between the surface normal $\n_p$ at $p$ and the surface normal $\n_q$ at any point $q$ on the (curved) edge $e\subset \S$ incident to $p$ in the given net of curvature lines on $\S$. To see this, consider an arc-length parametrized curve $\gamma:[0,\xi]\rightarrow \S$ with $\gamma(0)=p$ and $\gamma(\xi)=q$. Notice that our basic assumptions imply that $\gamma$ can be chosen such that $\xi \leq \epsilon$. The Gauss image $\tilde\gamma=\n\circ\gamma$ of $\gamma$ is a curve on the unit sphere. The length of $\tilde\gamma$ is therefore bounded from below by $\angle(\n_p,\n_q)$, \ie, the length of the minimizing geodesic joining $\n_p$ and $\n_q$ on $\mathbb{S}^2$. The tangent vector of $\tilde\gamma$ at a point $\tilde\gamma(s)$, $s\in [0,\xi]$, is given by $\W\gamma'(s)$, and the norm of this vector is therefore bounded above by $\K$. Hence,
\begin{align*}
\angle(\n_p,\n_q) \leq \int_0^\xi \| \W\gamma'(s) \| \,ds \leq \K\xi \leq \K\epsilon \leq \frac{1}{2} ~.
\end{align*}
Writing $q= (q_1, q_2, q_n)$ with respect to the Darboux frame centered at $p=(0,0,0)$, it follows that
\begin{align}\label{eqn:boundDh}
\|D h(q_1, q_2)\| = \tan \angle(\n_p,\n_q) \leq \tan\frac{1}{2} 
 ~.
\end{align}
Plugging this into~\eqref{eqn:d2vEstimate}, we find that $|D^2_\v h|\leq 2\K$, which, together with~\eqref{eqn:intH}, yields the first part of the lemma.

In order to prove the second part, we first note that $D^2h(0) = D^2 P(0)$ and that $D^2P$ is constant. Hence
\begin{align}\label{eqn:hPbound}
|h(\ex,\ey)-P(\ex,\ey)|=
\left|\int_0^d\int_0^t \int_0^\tau D_\v^3 h(\sigma \v) \, d\sigma \, d\tau \, dt \right|  \ .
\end{align}
We consequently seek a bound for $|D_\v^3 h|$ in terms of $\K$ and $\dK$. Considering~\eqref{eqn:kv} and taking another derivative with respect to $\v$ yields
\begin{align}
\begin{split}
\label{eqn:d3v}
D^3_\v h &= D_\v (\kappa_\v) \sqrt{1+\|Dh\|^2} \I(\v,\v) \\
& \quad + \kappa_\v D_\v \left(\sqrt{1+\|Dh\|^2}\right) \I(\v,\v) \\
& \quad + \kappa_\v \sqrt{1+\|Dh\|^2} D_\v (\I(\v,\v)) \ .
\end{split}
\end{align}
We bound the terms appearing on the right hand side one by one. To treat the first term, we use \eqref{eqn:kv} and $\nabla I\equiv0$, where $\nabla$ denotes covariant differentiation with respect to the metric induced by $\I$, to derive
\begin{align}
\label{eqn:dv_kv}
D_\v \kappa_\v &= \frac{\I((\nabla_\v \W)\v,\v) + 2 \I(\W\v,\nabla_\v \v) - 2 \kappa_\v \I(\v,\nabla_\v \v)}{\I(\v,\v)} ~.
\end{align}
Let $\u(x,y)=\u$ be a constant vector field in the parameter domain. Using the Koszul formula for the Levi-Civita connection and applying \eqref{eqn:I} yields
\begin{align*}
2 \I(\nabla_\v \v, \u)
= 2\v(\I(\v,\u)) - \u(\I(\v,\v))
= 2 D^2_\v h \, D_\u h \ .
\end{align*}
The last equality only depends on the value of $\u(x,y)$ at the point $(x,y)$ and therefore holds for any field $\u$ in the parameter domain. From~\eqref{eqn:d2vEstimate},~\eqref{eqn:boundDh}, and $\|\v\|=1$, we obtain
\begin{align*}
| \I(\nabla_\v \v, \u) |
\leq C \,  \|\u\| \ ,
\end{align*}
with $C= C(\K)$.
This can be used in \eqref{eqn:dv_kv}, together with our bounds on the norms of $\W$ and $\nabla\W$ in terms of $\K$ and $\dK$, respectively, to obtain an upper bound for the first term in~\eqref{eqn:d3v} by $C(\K,\dK)$. Bounds for the remaining two terms in~\eqref{eqn:d3v} can be obtained in a similar fashion using \eqref{eqn:I} and \eqref{eqn:kv}, proving the estimate $|D^3_\v h| \leq C(\K,\dK)$. Using~\eqref{eqn:hPbound} and $d\leq \epsilon$ then implies the claim of the second part in the statement of the lemma.
\end{proof}

\begin{corollary}
\label{edge normal component} With the same assumptions as in Lemma~\ref{lemma:3rdOrderBound} and by defining $\dk=\kappa_2-\kappa_1$, we obtain
\begin{align*}
\ez=\frac{\kx}{2}\e^2+ \frac{\dk}{2} \ey^2 + \O{\epsilon^3}
=\frac{\ky}{2}\e^2- \frac{\dk}{2} \ex^2 + \O{\epsilon^3} \ ,
\end{align*}
with $|\O{\epsilon^3}| \leq C(\K, \dK)\epsilon^3$.
\end{corollary}
\begin{proof}
From the second part of Lemma~\ref{lemma:3rdOrderBound} and simple algebraic manipulations, we deduce that
\begin{align*}
\ez= 
\frac{\kx}{2}\e^2 - \frac{\kx}{2} \ez^2 + \frac{\dk}{2} \ey^2
+ \O{\epsilon^3}
=
\frac{\kx}{2}\e^2 + \frac{\dk}{2} \ey^2 
+ \O{\epsilon^3} \ .
\end{align*}
The last equality follows from Lemma~\ref{lemma:3rdOrderBound} and the sampling condition 
$\epsilon\leq\frac{1}{2\K}$ (which is implied by our basic assumptions), which together ensure that
$\ez^2\leq\K^2\epsilon^4\leq\frac{\K}{2}\epsilon^3$. The second equation in the statement of the
corollary follows analogously.
\end{proof}


\begin{figure}[t]
\begin{center}
\includegraphics[scale=1.2]{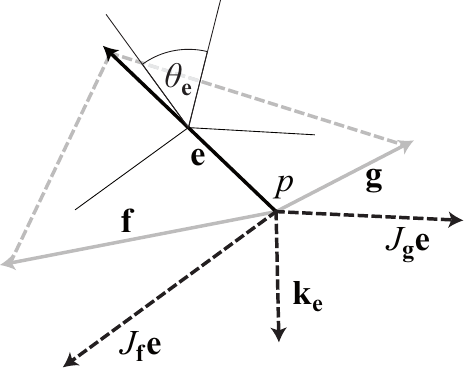}
\caption{The \emph{discrete curvature vector} $\ke$ is the sum of $\Jp\e$ (the rotation of $\e$ by $\frac\pi2$ about the axis $\e\times\f$) and $\Jm\e$ (the rotation of $\e$ by $\frac\pi2$ about the axis $\e\times\g$).}
\label{fig:curvatureVector}
\end{center}
\end{figure}

For the following discussion, it will be useful to work with the curvature \emph{vector} $\k$ corresponding to the discrete integrated curvature $\ke = 2\sin\frac{\te}{2} \|\e\|$ (see Section~\ref{subsec:dicrete.curvatures}). 
With our usual notions, consider the two edge vectors $\f$ and $\g$ directly neighboring $\e$ in $st(p)$ (see Figure~\ref{fig:curvatureVector}).
Recall from Corollary~\ref{cor:vertexStarProjection} that we may assume that $(\e\times\f) \cdot \n > 0$ and $(\e\times\g) \cdot \n < 0$, where $\n$ is the normal of $\S$ at $p$.
A straightforward
calculation reveals that we can express the discrete curvature vector by
\begin{align}\label{eqn:curvVector}
 \k = \Jp\e + \Jm\e \ ,%
\end{align}
where $\Jp$ and $\Jm$ denote rotations by $\frac\pi2$ around the axes
$\e\times\fp$ and $\e\times\fm$, respectively (see
Figure~\ref{fig:curvatureVector}). 
Consider now the splitting 
\begin{align*}
\k = (\k)_t+(\k)_n
\end{align*}
of $\k$ into its tangential and normal component with respect to $T_p\S$. We show that the tangential component is negligible:

\begin{lemma}\label{tangentialBound}
Under the same assumptions as in Lemma~\ref{lemma:3rdOrderBound}, the
projection of the discrete curvature vector $\k$ onto the tangent plane $T_p\S$ satisfies $\Vert(\k)_t\Vert\leq 10\K^2\rho^2\, \epsilon^3$.
\end{lemma}

\begin{proof}
With the notions and assumption of the preceding discussion, we first note that a straightforward calculation reveals that
\begin{align}\label{eqn:Je}
\J\e &= \frac{\e^2\f-(\f\cdot\e)\e}{\|\e\times\f\|} \ ,
\end{align}
and analogously for $\Jm\e$. 

With respect to $T_p\S$, let $\et$, $\ft$, and $\gt$ denote the tangential components of the vectors $\e$, $\f$, and $\g$, respectively. Then, by assumption, we have $(\et\times\ft) \cdot \n
> 0$ and $(\et\times\gt) \cdot \n < 0$. The tangential part of the discrete
curvature vector $\k$ is given by $(\k)_t = (\Jm\e)_t + (\Jp\f)_t$, where we deduce
from~\eqref{eqn:Je}, using the coordinates of our fixed Darboux frame at $p$, that
\begin{align*}
(\Jp\e)_t
= \frac{\|\e\|^2\ft-(\f\cdot\e)\et}{\|\e\times\f\|} 
= \frac{\|\et\times\ft\|}{\|\e\times\f\|}
(-\ey,\ex,0)+
\frac{\ez^2\ft-\ez\fz\et}{\|\e\times\f\|} \ ,
\end{align*}
and similarly
\begin{align*}
(\Jm\e)_t &= -\frac{\|\et\times\gt\|}{\|\e\times\g\|} (-\ey,\ex,0)+
\frac{\ez^2\gt-\ez\gz\et}{\|\e\times\g\|} \ .
\end{align*}
Using the first part of Lemma~\ref{lemma:3rdOrderBound} and the definition of shape regularity~\eqref{eqn:shapeRegularity}, we obtain
\begin{align*}
\|\e\times\f\|^2&=\|\et\times\ft\|^2+(\ez\fx-\ex\fz)^2+(\ez\fy-\ey\fz)^2 \\
&\leq\|\et\times\ft\|^2+8(\K\|\e\|\ \|\f\|\epsilon)^2 \\
&\leq\|\et\times\ft\|^2+8\K^2\rho^2\epsilon^2\|\e\times\f\|^2 \ .
\end{align*}
Therefore, we have
\begin{align*}
1 \geq \frac{\|\et\times\ft\|}{\|\e\times\f\|} \geq \frac{\|\et\times\ft\|^2}{\|\e\times\f\|^2}
\geq 1 - 8 \K^2\rho^2\epsilon^2 \ ,
\end{align*}
and similarly for $\g$.
It follows that the two terms in $(\Jp\e+\Jm\e)_t$ containing $(-\ey,\ex,0)$ cancel up to a term bounded by $(8\K^2\rho^2)\epsilon^3$, where the power of three is due to the fact that the norm of $(-\ey,\ex,0)$ is also bounded by $\epsilon$.
Moreover, we observe that the first part of Lemma~\ref{lemma:3rdOrderBound} and the definition of shape regularity~\eqref{eqn:shapeRegularity} yield
\begin{align*}
\frac{\|\ez^2\ft-\ez\fz\et\|}{\|\e\times\f\|}
\leq \frac{\K^2\|\e\|^4\|\f\| + \K^2\|\e\|^3\|\f\|^2}{\|\e\times\f\|}
\leq 2\K^2\rho\epsilon^3 \ ,
\end{align*}
and analogously for $\g$.
Therefore, we arrive at
\begin{align*}
\|(\k)_t\| =\|(\Jp\e+\Jm\e)_t\|\leq (10\K^2\rho^2) \epsilon^3 \ ,
\end{align*}
proving the claim.
\end{proof}

The above implies the main result of this section:

\begin{proof}[Proof of Proposition~\ref{lemma:edgeCurvatureBound}]
To estimate the normal component $(\k)_n=(\Jp\e)_n+(\Jm\e)_n$ of
the discrete curvature vector, we use~\eqref{eqn:Je} to obtain
\begin{align*}
(\J\e)_n &= \frac{\e^2\fz-(\f\cdot\e)\ez}{\|\e\times\f\|} \ ,
\end{align*}
and analogously for $(\Jm\e)_n$.
We note that the circumcentric area of $\e$ can be expressed in a similar manner as a sum $A_\e = A_{\e,\fp} +
A_{\e,\fm}$, where
\begin{align*}
A_{\e,\f}=\frac{\f\cdot(\f-\e)}{4\|\e\times\f\|}\e^2 \ ,
\end{align*}
and analogously for $A_{\e,\g}$. Applying Corollary~\ref{edge normal component}, 
we obtain
\begin{align*}
(\Jp\e)_n
&= \frac{\e^2(\kappa_2\f^2-\dk \fx^2)-\f\cdot\e(\ky\e^2-\dk \ex^2)}{2\|\e\times\f\|} + \O{\epsilon^3}\\
&= 2A_{\e,\f} \ky
-\dk\fx^2\frac{\e^2  }{2\|\e\times\f\|}
-\dk\f\cdot\e\frac{\ex^2}{2\|\e\times\f\|}
+\O{\epsilon^3} ~,
\end{align*}
and analogously for $(\Jm\e)_n$, with $|\O{\epsilon^3}|\leq C(\K, \dK, \rho)\epsilon^3$.  Applying the shape regularity condition
(\ref{eqn:shapeRegularity}) yields
\begin{align*}
\left| (\J\e)_n - 2A_{\e,\f} \ky \right| &\leq C
\left(\left|\dk(\fx^2+\f\cdot\e) \right| +\epsilon^3 \right) \ ,
\end{align*}
and similarly for the difference between $(\Jm\e)_n$ and $2A_{\e,\fm}\ky$.
According to the statement of Proposition~\ref{lemma:edgeCurvatureBound}, 
we assume that $\e$ is associated with the direction $\vx$, whereas $\f$ and $\g$ are associated with $\vy$. Using our notion of \emph{tangential deviation} from Definition~\ref{def:tangentialDeviation} and noting that $| \f\cdot\e| 
\leq | \ey+\fx |\epsilon +\O{\epsilon^4}$, we arrive at
\begin{align*}
\left| (\k)_n - 2A_\e \ky \right| \leq C (\left| \dk(\ed+\fd+\gd) \right|
\epsilon + \epsilon^3) \ ,
\end{align*}
which implies the claim, since $(\k)_t$ is of order
$\O{\epsilon^3}$ by Lemma~\ref{tangentialBound}.
\end{proof}
%
%
%
%

\subsection{Estimates for non-umbilical vertices}
\label{subsec:deltaKappa} 

As mentioned before, the results of the preceding section are not entirely specific for nets of curvature lines but hold for a larger class of smoothly embedded nets. As such, these results do not suffice to prove our main error estimate in Theorem~\ref{thm:main}, due to the failure of uniform convergence of curvature approximations constructed in a $1$-local manner for general nets. Indeed, assuming that $\e$ is associated with the principal direction $\vx$, we seek a bound of the form
\begin{align}\label{eqn:mainEstimate}
\left| \ke - 2A_\e \ky \right|\leq C\epsilon^3 \ ,
\end{align}
from which we may derive the desired main error estimate by employing our results on the existence of uniform lower and upper bounds of (sufficiently many) circumcentric areas (see Section~\ref{subsec:area}). 

While a bound of the form~\eqref{eqn:mainEstimate} does not hold for general smooth nets, it is indeed valid for nets of curvature lines. This is a consequence of the results of the preceding section and the fact that for nets of curvature lines we have
\begin{align}\label{eqn:deltaKappa} 
|\dk(p)\ed|\leq C\epsilon^2 \ ,
\end{align}
which is trivially satisfied for umbilical vertices and true for non-umbilical ones \emph{provided} that we associate $\e$ with its canonical principal direction. This is precisely the main result of this section.

Observe that there is a simple case, where~\eqref{eqn:deltaKappa} is obviously fulfilled:
let $\S$ be a paraboloid, and let $p$ be its apex. Consider the four edge vectors emanating from $p$ in a local polyhedral approximation of a net of curvature lines containing $p$ as a vertex. Then the tangential deviation of each of these edges vanishes, so~\eqref{eqn:deltaKappa} is clearly satisfied.

The main difference between nets of curvature
lines on arbitrary smooth surfaces and the specific case of a paraboloid is the
fact that curvature lines usually have non-zero geodesic curvature $\kappa^g$. While
the tangential deviation $\ed$ can always be bounded by $C\epsilon^2$, with $C= C(\K,\kappa^g)$, this does not suffice for a uniform error
bound, since $\kappa^g$ may blow up at umbilical points. Perhaps surprisingly,
though, the product $|\dk(p)\ed|$ {\em can} be bounded for nets of curvature
lines: 
%

\begin{proposition}\label{lemma:tangentDeviation}
Let $\e$ be an edge vector emanating from a non-umbilical vertex $p$ in a local polyhedral approximation of a discrete net of curvature lines on $\S$.
If $\e$ is associated with its canonical principal direction, then our basic assumptions imply
\begin{align*}
|\dk(p)\ed|  \leq (\K^2+4\dK){\epsilon^2} \ ,
\end{align*}
where $\ed$ is the tangential deviation of $\e$ and $\K$, $\dK$ denote our usual bounds on normal curvatures and their derivatives, respectively.
\end{proposition}
The intuition behind this statement is as follows. Roughly, $\ed$ is
proportional to the geodesic curvature of the curvature line corresponding to
$\e$. This geodesic curvature, in turn, is inversely proportional to $\dk$, as
the next lemma shows. Therefore, the product $\dk\ed$ can be uniformly
bounded.


\begin{lemma}\label{lemma:geodCurv}
At any non-umbilical point of $\S$, the geodesic curvature $\kgx$ of the principal curvature line along $\vx$ satisfies
\begin{align*}
\kgx = \frac{\nabla_{\vy}\kx} {\kx-\ky} \quad\text{and thus}\quad%
|\kgx|\leq \frac{\dK} {|\kx-\ky|} \ ,
\end{align*}
where $\kx$ and $\ky$ denote the principal curvatures corresponding to the
principal directions $\vx$ and $\vy$, respectively.
\end{lemma}
\begin{proof}
Since $\vx$ and $\vy$ are orthonormal eigenvectors of the shape operator
$\W$, we have $\W\vx\cdot\vy=0$. We use the Frenet formulas
$\nabla_{\vx}\vx=\kgx\vy$, $\nabla_{\vx}\vy=-\kgx\vx$, $\nabla_{\vy}\vx=\kgy\vy$, and $\nabla_{\vy}\vy=-\kgy\vx$ as
well as the Codazzi-Mainardi equation $(\nabla_\u\W)\v=(\nabla_\v\W)\u$ to
obtain
\begin{equation*}
\begin{split}
0&=
\nabla_{\vx}(\W\vx\cdot\vy)=
\nabla_{\vx}(\W\vy\cdot\vx)
\\&=
(\nabla_{\vx}\W)\vy\cdot\vx  +
\W(\nabla_{\vx}\vy)\cdot\vx +
\W\vy\cdot\nabla_{\vx}\vx
\\&=
(\nabla_{\vy}\W)\vx\cdot\vx +
\W(-\kgx\vx)\cdot\vx +
\W\vy\cdot(\kgx\vy)
\\&=
(\nabla_{\vy}\W)\vx\cdot\vx +
\kgx(-\kappa_1\vx\cdot\vx +\kappa_2\vy\cdot\vy)
\\&=
\nabla_{\vy} (\W\vx\cdot\vx) - \W(\nabla_{\vy}\vx)\cdot\vx -
\W\vx\cdot\nabla_{\vy}\vx + \kgx(\ky-\kx)
\\&=
\nabla_{\vy} (\W\vx\cdot\vx) - \W(\kgy\vy)\cdot\vx -
\W\vx\cdot\kgy\vy + \kgx(\ky-\kx)
\\&=
\nabla_{\vy}\kx + \kgx(\ky-\kx) \ ,
\end{split}
\end{equation*}
proving the first part. The second part follows from the definition of $\dK$.
\end{proof}


\begin{proof}[Proof of Proposition~\ref{lemma:tangentDeviation}]
Let $\gamma:[0,\epsilon]\to \S$ be the curvature line that is canonically associated with $\e$,
parameterized by arc-length and passing through $p=\gamma(0)$. By definition,
$\ed$ is bounded above by the maximum distance from $\gamma$ to the tangent line
passing trough $\gamma'(0)$. Since $\gamma$ is parameterized by arc-length, we
obtain
\begin{align}\label{eqn:edKappa}
    \ed \leq \frac{\K_\gamma}{2} \epsilon^2 \ ,
\end{align}
where $\K_\gamma$ denotes the maximum curvature of $\gamma$ as a space curve.
Decomposing the curvature vector of $\gamma$ into its normal and geodesic 
components, and denoting by $\K^n_\gamma$ and $\K^g_\gamma$ the respective maxima of the norms of these components,
Lemma~\ref{lemma:geodCurv} yields
\begin{align}\label{eqn:KGamma}
\K_\gamma \leq \K^n_\gamma + \K^g_\gamma \leq \K + \max_{s\in [0,
\epsilon]}\frac{\dK}{|\dk(\gamma(s))|}\ ,
\end{align}
since $\K^n_\gamma \leq \K$ by definition of $\K$.  Consequently, we seek a lower bound for
$|\dk(\gamma(s))|$. To do so, first observe that by definition of $\dK$ the derivatives of $\kx$ and $\ky$ are
bounded by $\dK$. Hence, the function $\dk$ is Lipschitz with constant
$2\dK$, \ie,
\begin{align*}
|\dk(p)-\dk(q)|\leq 2\dK  d_\S(p,q) \ ,
\end{align*}
for every point $q\in\S$.
We now distinguish two cases: (i) $|\dk(p)| < 4\dK \epsilon$ and
(ii)~$|\dk(p)|\geq 4\dK \epsilon$. In the first case, we immediately obtain
\begin{align*}
|\dk(p) \ed| \leq 4\dK \epsilon^2 \ ,
\end{align*}
which already proves the claim of the lemma. In the second case, we observe that for all $s\in [0,\epsilon]$, we have
\begin{align*}
|\dk(\gamma(s))| \geq |\dk(p)|-2\dK \epsilon \geq \frac{1}{2}|\dk(p)| \ .
\end{align*}
Plugging this into~\eqref{eqn:KGamma} gives
\begin{align*}
\K_\gamma \leq \K + \frac{2\dK}{|\dk(p)|} \ .
\end{align*}
Together with~\eqref{eqn:edKappa}, this yields
\begin{align*}
|\dk(p) \ed| &\leq |\dk(p)|
\left(\frac{\K}{2}+\frac{\dK}{|\dk(p)|}\right)\epsilon^2 \leq
(\K^2+\dK)\epsilon^2 ~,
\end{align*}
completing the proof.
\end{proof}

\subsection{Combining the strings}\label{sec:wrapUp}
We are now in the position to prove Theorem~\ref{thm:main}. 

To this end, assume that the edge $\e$ in our local polyhedral approximation is associated with the principal direction given by $\vx$, and let this be the canonical direction if $p$ is not umbilical. Furthermore, let $\ke = 2\sin\frac{\te}{2} \|\e\|$ be the discrete integrated curvature introduced in Section~\ref{subsec:dicrete.curvatures}. Then propositions~\ref{lemma:edgeCurvatureBound} and~\ref{lemma:tangentDeviation} imply
that there exists a constant $C=C(\K,\dK,\rho)$ such that
\begin{align}\label{eqn:discreteCurvEstimate}
\left| \ke - 2\Ae \ky \right| \leq C\epsilon^3 \ .
\end{align}
Assume additionally that $\e$ is chosen such that it maximizes the circumcentric area $\Ae$. (There are exactly two choices for non-umbilical vertices.) Proposition~\ref{cor:largestEdgeIsFine} shows that this choice leads to lower and upper bounds for $\Ae$ by $(1/C) \epsilon^2$ and $C \epsilon^2$, respectively, with
$C= C(\rho)$. Therefore, we can divide~\eqref{eqn:discreteCurvEstimate} by
$2\Ae$ to obtain
\begin{align*} \left|\frac{\ke}{2\Ae}-\ky\right| \leq C \epsilon \ ,
\end{align*}
with $C=C(\K,\dK,\rho)$. 

Finally, in order to prove our error estimate for the two other edge-based integrated discrete curvatures,  $\ke = 2\tan\frac{\te}{2} \|\e\|$ and $\ke = \te\|\e\|$, we infer from Lemma~\ref{lemma:normalBound} and a simple application of Taylor's theorem that there exists a constant $C=C(\K,\rho)$ such that
\begin{align*}
\left|\te - 2\sin\frac{\te}{2}\right| \leq C\epsilon^3\quad\text{and}\quad\left|\te - 2\tan\frac{\te}{2}\right| \leq C\epsilon^3 \ ,
\end{align*}
which can be used to obtain the requisite bound~\eqref{eqn:discreteCurvEstimate} for the other two discrete curvature definitions as well.
This completes the proof of Theorem~\ref{thm:main}. 

\paragraph*{Acknowledgments}
We would like to thank the anonymous reviewers for their very valuable and detailed comments. We also thank Emanuel Huhnen-Venedey und Ramsay Dyer for their helpful feedback. This work was partially supported by the DFG Research Center \textsc{Matheon} and the DFG Research Unit Polyhedral Surfaces.

\bibliography{CurvatureLineCotan}

\begin{thebibliography}{25}
\providecommand{\natexlab}[1]{#1}
\providecommand{\url}[1]{\texttt{#1}}
\expandafter\ifx\csname urlstyle\endcsname\relax
  \providecommand{\doi}[1]{doi: #1}\else
  \providecommand{\doi}{doi: \begingroup \urlstyle{rm}\Url}\fi

\bibitem[Belkin et~al.(2008)Belkin, Sun, and Wang]{belkin08laplace}
M.~Belkin, J.~Sun, and Y.~Wang.
\newblock \href{http://dx.doi.org/10.1145/1377676.1377725}{Discrete {L}aplace
  operator on meshed surfaces}.
\newblock In \emph{Proceedings SCG '08}, pages 278--287, 2008.

\bibitem[Berry and Hannay(1977)]{berry77umbilic}
M.~V. Berry and J.~H. Hannay.
\newblock \href{http://dx.doi.org/10.1088/0305-4470/10/11/009}{Umbilic points
  on {G}aussian random surfaces}.
\newblock \emph{J. Phys. A: Math. Gen}, 10:\penalty0 1809--1821, 1977.

\bibitem[Bobenko and Springborn(2007)]{BobenkoSpringborn05}
A.~Bobenko and B.~Springborn.
\newblock \href{http://dx.doi.org/10.1007/s00454-007-9006-1}{A discrete
  {L}aplace–{B}eltrami operator for simplicial surfaces}.
\newblock \emph{Discr.~Comp.~Geom.}, 38\penalty0 (4):\penalty0 740--756,
  December 2007.

\bibitem[Bobenko et~al.(2008)Bobenko, Schr{\"o}der, Sullivan, and
  Ziegler]{DDGBook}
A.~Bobenko, P.~Schr{\"o}der, J.~Sullivan, and G.~Ziegler, editors.
\newblock \emph{{Discrete Differential Geometry}}, volume~38 of
  \emph{Oberwolfach Seminars}.
\newblock Birkh{\"a}user, 2008.

\bibitem[Bobenko and Suris(1999)]{bobenko99lagrange}
A.~I. Bobenko and Y.~Suris.
\newblock \href{http://dx.doi.org/10.1007/s002200050642}{Discrete time
  {L}agrangian mechanics on {L}ie groups, with an application to the {L}agrange
  top}.
\newblock \emph{Comm. Math. Phys.}, 204\penalty0 (1):\penalty0 147--188, 1999.

\bibitem[Bobenko and Suris(2008)]{Bobenko2008Discrete}
A.~I. Bobenko and Y.~B. Suris.
\newblock \emph{Discrete Differential Geometry: Integrable Structure},
  volume~98 of \emph{Graduate Studies in Mathematics}.
\newblock American Mathematical Society, 2008.

\bibitem[Cazals and Pouget(2005)]{CazalsPouget05}
F.~Cazals and M.~Pouget.
\newblock \href{http://dx.doi.org/10.1016/j.cagd.2004.09.004}{Estimating
  differential quantities using polynomial fitting of osculating jets}.
\newblock \emph{Comput. Aided Geom. Des.}, 22\penalty0 (2):\penalty0 121--146,
  2005.

\bibitem[Cohen-Steiner and Morvan(2006)]{MorvanCohen-Steiner06}
D.~Cohen-Steiner and J.-M. Morvan.
\newblock
  \href{http://www.intlpress.com/JDG/p/2006/74_3/JDG-74-3-363-394.pdf}{Second
  fundamental measure of geometric sets and local approximation of curvatures}.
\newblock \emph{J. Differential Geom.}, 73\penalty0 (3):\penalty0 363--394,
  2006.

\bibitem[Desbrun et~al.(2005)Desbrun, Hirani, Leok, and Marsden]{desbrun05dec}
M.~Desbrun, A.~N. Hirani, M.~Leok, and J.~E. Marsden.
\newblock \href{http://arxiv.org/abs/math/0508341}{Discrete exterior calculus}.
\newblock Arxiv preprint, 2005.
\newblock \href {http://arxiv.org/abs/math/0508341}
  {\path{arXiv:math/0508341}}.

\bibitem[Dyer and Schaefer(2009)]{Dyer2009Circumcentric}
R.~Dyer and S.~Schaefer.
\newblock \href{ftp://fas.sfu.ca/pub/cs/TR/2009/CMPT2009-06.pdf}{Circumcentric
  dual cells with negative area}.
\newblock Technical Report TR 2009-02, School of Computing Science, Simon
  Fraser University, Burnaby, BC, Canada, 2009.

\bibitem[Dziuk(1988)]{dziuk88laplace}
G.~Dziuk.
\newblock \href{http://dx.doi.org/10.1007/BFb0082859}{Finite elements for the
  {B}eltrami operator on arbitrary surfaces}.
\newblock In \emph{Partial differential equations and calculus of variations},
  volume 1357 of \emph{Lecture Notes in Mathematics}, pages 142--155. Springer,
  1988.

\bibitem[Federer(1959)]{Federer59}
H.~Federer.
\newblock \href{http://www.jstor.org/stable/1993504}{Curvature measures}.
\newblock \emph{Trans. Amer. Math.}, 93\penalty0 (3):\penalty0 418--491, 1959.

\bibitem[Fu(1993)]{fu93convergence}
J.~H.~G. Fu.
\newblock
  \href{http://www.intlpress.com/JDG/archive/1993/37-1-177.pdf}{{Convergence of
  curvatures in secant approximations}}.
\newblock \emph{J. Differential Geom.}, 37:\penalty0 177--190, 1993.

\bibitem[Grinspun and Desbrun(2006)]{SiggraphCourseNotes}
E.~Grinspun and M.~Desbrun, editors.
\newblock \href{http://dx.doi.org/10.1145/1185657}{\emph{{Discrete differential
  geometry: an applied introduction}}}, ACM SIGGRAPH Courses Notes, 2006. ACM
  Press New York, NY, USA.

\bibitem[Hildebrandt et~al.(2006)Hildebrandt, Polthier, and
  Wardetzky]{wardetzky06convergence}
K.~Hildebrandt, K.~Polthier, and M.~Wardetzky.
\newblock \href{http://dx.doi.org/10.1007/s10711-006-9109-5}{{On the
  convergence of metric and geometric properties of polyhedral surfaces}}.
\newblock \emph{Geom. Dedicata}, 123\penalty0 (1):\penalty0 89--112, 2006.

\bibitem[Hoffmann(2008)]{hoffmann08hashimoto}
T.~Hoffmann.
\newblock \href{http://arxiv.org/abs/math/0007150}{Discrete {H}ashimoto
  surfaces and a doubly discrete smoke-ring flow}.
\newblock In \emph{Discrete Differential Geometry}, volume~38 of
  \emph{Oberwolfach Seminars}, pages 95--116. Birkh{\"a}user, 2008.

\bibitem[Meek and Walton(2000)]{MeekWalton00}
D.~S. Meek and D.~J. Walton.
\newblock \href{http://dx.doi.org/10.1016/S0167-8396(00)00006-6}{On surface
  normal and {G}aussian curvature approximations given data sampled from a
  smooth surface.}
\newblock \emph{Comput. Aided Geom. Des.}, 17\penalty0 (6):\penalty0 521--543,
  2000.

\bibitem[Morvan and Thibert(2004)]{MorvanThibert04}
J.-M. Morvan and B.~Thibert.
\newblock \href{http://dx.doi.org/10.1007/s00454-004-1096-4}{Approximation of
  the normal vector field and the area of a smooth surface}.
\newblock \emph{Discr.~Comp.~Geom.}, 32\penalty0 (3):\penalty0 383--400, 2004.

\bibitem[Pinkall and Polthier(1993)]{pinkall93minimal}
U.~Pinkall and K.~Polthier.
\newblock \href{http://projecteuclid.org/euclid.em/1062620735}{Computing
  discrete minimal surfaces and their conjugates}.
\newblock \emph{Exp. Math.}, 2\penalty0 (1):\penalty0 15--36, 1993.

\bibitem[Steiner(1840)]{steiner40parallel}
J.~Steiner.
\newblock {Ueber parallele Fl\"{a}chen.}
\newblock \emph{Monatsbericht der Akademie der Wissenschaften zu Berlin}, pages
  114--118, 1840.

\bibitem[Sullivan(2008)]{sullivan08curves}
J.~M. Sullivan.
\newblock \href{http://arxiv.org/abs/math/0606007}{Curves of finite total
  curvature}.
\newblock In \emph{Discrete Differential Geometry}, volume~38 of
  \emph{Oberwolfach Seminars}, pages 137--162. Birkh{\"a}user, 2008.

\bibitem[Wintgen(1982)]{Wintgen82}
P.~Wintgen.
\newblock Normal cycle and integral curvature for polyhedra in {R}iemannian
  manifolds.
\newblock In So\'{o}s and Szenthe, editors, \emph{Differential Geometry}, pages
  805--816. North-Holland, Amsterdam, 1982.

\bibitem[Xu(2006)]{xu06sphere}
G.~Xu.
\newblock \href{http://dx.doi.org/10.1142/S0218195906001938}{Discrete
  {L}aplace-{B}eltrami operator on sphere and optimal spherical
  triangulations}.
\newblock \emph{Int. J. Comp. Geom. Appl.}, 16\penalty0 (1):\penalty0 75--93,
  2006.

\bibitem[Xu et~al.(2005)Xu, Xu, and Sun]{xu05convergence}
Z.~Xu, G.~Xu, and J.-G. Sun.
\newblock \href{http://dx.doi.org/10.1007/11537908_27}{{Convergence analysis of
  discrete differential geometry operators over surfaces}}.
\newblock In \emph{Mathematics of Surfaces XI}, volume 3604 of \emph{LNCS},
  pages 448--457. Springer, 2005.

\bibitem[Z\"{a}hle(1986)]{Zahle86}
M.~Z\"{a}hle.
\newblock \href{http://dx.doi.org/10.1007/BF01195026}{Integral and current
  representations of {F}ederer's curvature measures}.
\newblock \emph{Arch. Math.}, 46:\penalty0 557--567, 1986.

\end{thebibliography}

\end{document}